\documentclass[12pt,a4paper,russian,titlepage]{article}
\usepackage{babel,amsmath,amssymb,amsthm,graphicx,amsfonts,latexsym,epsfig,longtable,wrapfig}

\newtheorem{theorem}{Теорема}
\newtheorem{lemma}{Лемма}

\newcommand{\EHR}{{\mathrm{EHR}}}

\begin{document}

\begin{center}{\Large Расширение $k$-закона нуля или единицы}
\end{center}

\vspace{0.1cm}

\begin{center}{\large М.Е. Жуковский}
\end{center}

\vspace{0.5cm}

\section{Введение}
\label{intro}

В данной работе изучаются предельные вероятности свойств первого
порядка случайного графа в модели Эрдеша--Реньи. Ранее нами был
установлен $k$-закон нуля или единицы, который описывает поведение
вероятностей свойств первого порядка, выраженных формулами с
ограниченной числом $k$ кванторной глубиной (см. \cite{Zhuk2}).
Закон выполнен для случайного графа $G(N,N^{-\alpha})$ при
$\alpha\in(0,1/(k-2))$. В настоящей работе мы расширили упомянутый
результат, увеличив диапазон значений $\alpha$, при которых
справедлив $k$-закон. В этом разделе мы сформулируем основной
результат, но прежде дадим все необходимые определения, а также
опишем кратко историю исследуемой задачи.

Пусть $N\in\mathbb{N}$, $0\leq p\leq 1.$ Рассмотрим множество
$\Omega_N=\{G=(V_N,E)\}$ всех неориентированных графов без петель
и кратных ребер с множеством вершин $V_N=\{1,...,N\}$.
\emph{Случайный граф в модели Эрдеша--Реньи} (см.
\cite{Erdos}--\cite{Alon}) это случайный элемент $G(N,p)$ со
значениями во множестве $\Omega_N$ и распределением ${\sf
P}_{N,p}$ на $\mathcal{F}_N=2^{\Omega_N},$ определенным формулой
$$
 {\sf P}_{N,p}(G)=p^{|E|}(1-p)^{C_N^2-|E|}.
$$

Говорят, что случайный граф {\it подчиняется закону нуля или
единицы} для класса свойств $\mathcal{C}$, если вероятность
выполнения каждого свойства из этого класса стремится либо к $0$,
либо к $1$.

Пожалуй, самым изученным в этом направлении является класс
свойств, выражаемых формулами первого порядка (см. \cite{Veresh},
\cite{Usp}). Эти формулы строятся с помощью символов отношения
$\sim,=$; логических связок
$\neg,\Rightarrow,\Leftrightarrow,\vee,\wedge$; переменных
$x,y,x_1...$; кванторов $\forall,\exists$. Символы переменных в
языке первого порядка обозначают вершины графа. Символ отношения
$\sim$ выражает свойство двух вершин быть соединенными ребром,
символ отношения $=$ выражает свойство двух вершин совпадать.
Обозначим $\mathcal{P}$ множество функций $p=p(N)$, для которых
выполнен закон нуля или единицы для класса свойств первого
порядка. Класс свойств первого порядка обозначим $\mathcal{L}.$ В
1969 году Ю.В.~Глебский, Д.И.~Коган, М.И.~Лиогонький и
В.А.~Таланов  (см. \cite{Kogan}) и независимо в 1976 году Р.~Фагин
(см. \cite{Fagin}) доказали, что если
$$
 \forall \alpha>0 \,\,\,\, N^{\alpha}\min\{p,1-p\}\rightarrow\infty,
 \,\, N\rightarrow\infty,
$$
то $p\in\mathcal{P}$. Кроме того, в 1988 году Дж. Спенсер и С.
Шела (см. \cite{Shelah}, \cite{Strange}) установили, что все
функции $p=N^{-\alpha}$, $\alpha\in\mathbb{R}\setminus\mathbb{Q}$,
$\alpha\in(0,1)$, также принадлежат $\mathcal{P}$. Разумеется,
$p=1-N^{-\alpha}\in\mathcal{P}$, если
$\alpha\in\mathbb{R}\setminus\mathbb{Q}$, $\alpha\in(0,1)$.

Заметим, что если $\alpha$
--- рациональное, $0<\alpha\leq 1$ и $p=N^{-\alpha}$, то несложно
показать, что случайный граф $G(N,p)$ не подчиняется закону нуля
или единицы (см. \cite{Alon}).

Помимо класса $\mathcal{L}$ рассматривается класс свойств
$\mathcal{L}^{\infty}\supset\mathcal{L}$, выражаемых формулами,
которые могут содержать бесконечное количество конъюнкций или
дизъюнкций (предложения в таком языке, вообще говоря, бесконечны,
в отличие от предложений, выражающих свойства из класса
$\mathcal{L}$). В 1997 году М.~Мак\-Артур (см. \cite{McArtur})
были установлены некоторые законы нуля или единицы для класса
свойств $\mathcal{L}^{\infty}_k\subset \mathcal{L}^{\infty}$,
выражаемых формулами с кванторной глубиной, ограниченной числом
$k$. А именно, были получены законы для случайного графа $G(N,p)$
при $p=N^{-\alpha}$ и некоторых рациональных $\alpha$ из $(0,1]$.
Рассмотрим теперь класс $\mathcal{L}_k\subset \mathcal{L}$
свойств, выражаемых формулами с кванторной глубиной, ограниченной
числом $k$ (предложения конечны, в отличие от предложений,
задающих свойства из класса $\mathcal{L}_k^{\infty}$). В 2012 году
мы доказали, что при $k\geq 3$ и $\alpha\in(0,1/(k-2))$ (см.
\cite{Zhuk2}, \cite{zhuk_dan}) случайный граф $G(N,N^{-\alpha})$
подчиняется закону нуля или единицы для класса $\mathcal{L}_k$
(если случайный граф подчиняется закону нуля или единицы для
класса $\mathcal{L}_k$, то мы говорим, что он {\it подчиняется
$k$-закону нуля или единицы}). Кроме того, в упомянутой работе мы
доказали, что при $\alpha=1/(k-2)$ случайный граф
$G(N,N^{-\alpha})$ не подчиняется $k$-закону нуля или единицы. В
\cite{zhuk_MJC} установлено, что существуют пределы вероятностей
${\sf P}_{N,N^{-1/(k-2)}}$ выполнения всех свойств, выражаемых
формулами
с кванторной глубиной, ограниченной числом $k$.\\

Сформулируем основные результаты данной работы. Далее мы
рассматриваем только случай $k>3$, так как при $k=3$ (а,
следовательно, и при $k<3$) случайный граф $G(N,N^{-\alpha})$
подчиняется $k$-закону нуля или единицы для всех $\alpha\in(0,1)$
(см. \cite{zhuk_MJC}).

\begin{theorem}

Пусть $k>3$ --- произвольное натуральное число. Пусть, кроме того,
$\mathcal{Q}$ --- множество положительных дробей с числителем, не
превосходящим числа $2^{k-1}$. Случайный граф $G(N,N^{-{\alpha}})$
подчиняется $k$-закону нуля или единицы, если
$\alpha=1-\frac{1}{2^{k-1}+\beta}$, $\beta\in(0,\infty)\setminus
\mathcal{Q}$.

\label{positive}
\end{theorem}

Итак, мы рассмотрели интервал $(1-2^{1-k},1)$ и получили множество
рациональных чисел $\alpha$ из этого интервала, при которых
$k$-закон справедлив. Так как любое число из $(1-2^{1-k},1)$
представляется в виде $1-\frac{1}{2^{k-1}+\beta}$, закон будет
выполнен при любых $\alpha$ из
$$
 \left(1-\frac{1}{2^k},1\right)
 \bigcup\left(1-\frac{1}{2^k-1},1-\frac{1}{2^k}\right)
 \bigcup\ldots\bigcup\left(1-\frac{1}{2^{k-1}+
 2^{k-2}},1-\frac{1}{2^{k-1}+2^{k-2}+1}\right)
 \bigcup
$$
$$
 \left(1-\frac{1}{2^{k-1}+\frac{2^{k-1}-1}{2}},1-\frac{1}{2^{k-1}+2^{k-2}}\right)
 \bigcup\ldots
 \bigcup
$$
$$
 \left(1-\frac{1}{2^{k-1}+\frac{2^{k-1}}{3}},1-\frac{1}{2^{k-1}+\frac{2^{k-1}-\left[\frac{2^{k-1}}{3}\right]}{2}}\right)\bigcup\ldots.
$$
Длины интервалов уменьшаются при стремлении концов к $1-2^{1-k}$.
Кроме того, для некоторых из них мы доказали, что на концах
интервалов $k$-закон не выполнен (см.
теорему \ref{negative}).\\

Доказательство теоремы будет построено по следующей схеме: в
разделе \ref{game} мы сформулируем теорему Эренфойхта, связывающую
законы нуля или единицы и существование выигрышной стратегии
второго игрока в игре Эренфойхта. В том же разделе мы определим
стратегию второго игрока, которая с вероятностью, стремящейся к 1,
приведет его к победе, как мы покажем в разделе \ref{proof1}.
Раздел \ref{constructions} посвящен построению конструкций,
свойства которых позволяют нам доказать, что второй игрок "почти
всегда"$\,\,$сможет играть в соответствии с
определенной нами стратегией.\\

Сформулируем результат о нарушении закона при
$\alpha=1-\frac{1}{2^{k-1}+\beta}$ для некоторых $\beta$ из
множества $\mathcal{Q}$, определенного в  теореме \ref{positive}.

\begin{theorem} Пусть $k>3$ --- произвольное натуральное число. Пусть, кроме того,
$\tilde{\mathcal{Q}}$ --- множество натуральных чисел, не
превосходящих $2^{k-1}$. Случайный граф $G(N,N^{-{\alpha}})$ не
подчиняется $k$-закону нуля или единицы, если
$\alpha=1-\frac{1}{2^{k-1}+\beta}$, $\beta\in
\tilde{\mathcal{Q}}$.

\label{negative}
\end{theorem}

Доказательство теоремы мы приводим в разделе \ref{proof2}.
Подобные рассуждения мы можем проводить и для некоторых дробных
$\beta$, но общего результата наш метод не дает, поэтому в данной
работе рассматриваются только целые $\beta$.






\section{Конструкции}
\label{constructions}

Данный раздел мы начнем с описания известных результатов о
распределении малых подграфов в случайном графе Эрдеша--Реньи. Для
простоты изложения мы приводим в этом разделе все утверждения,
которые будут использованы для доказательства основных результатов
работы. В параграфе \ref{cicles} будет сформулирована и доказана
лемма, играющая ключевую роль при доказательстве теоремы
\ref{positive}.

\subsection{Малые подграфы и расширения}
\label{subgraphs}

Для произвольного графа $G$ обозначим $v(G)$ число его вершин,
$e(G)$ --- число его ребер, $a(G)$ --- число его автоморфизмов.
Плотность графа $\frac{e(G)}{v(G)}$ обозначим $\rho(G)$. Граф $G$
называется \emph{сбалансированным}, если для каждого его подграфа
$H$ выполнено неравенство $\rho(H)\leq \rho(G).$ Граф $G$
\emph{строго сбалансированный}, если для любого $H\subset G$
справедливо строгое неравенство $\rho(H)<\rho(G).$ Сформулируем
теорему (см. \cite{Erdos}--\cite{Bol}, \cite{Alon}) о количестве
копий строго сбалансированного графа. Пусть $N_G$ --- количество
копий $G$ в случайном графе $G(N,p)$. Положим
$\rho^{\max}(G)=\max\{\rho(H):\,
H\subseteq G\}.$\\

\begin{theorem} [\cite{Erdos}] Пусть $G$ --- произвольный граф. Если $p=o\left(N^{-1/\rho^{\max}(G)}\right)$,
то
$$
 \lim_{N\rightarrow\infty}{\sf P}_{N,p}(N_G>0)=0.
$$
Если же
$N^{-1/\rho^{\max}(G)}=o(p)$, то
$$
 \lim_{N\rightarrow\infty}{\sf P}_{N,p}(N_G>0)=1.
$$
Пусть теперь $G$
--- строго сбалансированный граф. Если $N^{-1/\rho(G)}=o(p),$ то
для любого $\varepsilon>0$ справедливо равенство
$$
\lim\limits_{N\rightarrow\infty}{\sf P}_{N,p}(\left|N_G- {\sf
E}_{N,p}N_G\right|\leq \varepsilon {\sf E}_{N,p}N_G)=1,
$$
где ${\sf E}_{N,p}$ --- это математическое ожидание по мере ${\sf
P}_{N,p}$ (здесь и далее в подобных ситуациях мы будем говорить,
что с асимптотической вероятностью $1$ выполнено $N_G\sim {\sf
E}_{N,p}N_G$). Если же $p=N^{-1/\rho(G)},$ то
$$
\lim\limits_{N\rightarrow\infty}{\sf P}_{N,p}(N_G=0)=e^{-1/a(G)}.
$$
\label{erdos}
\end{theorem}

Обратимся к задаче, поставленной Дж.~Спенсером в 1990 году (см.
\cite{Alon}, \cite{Spencer}). Рассмотрим такие графы
$H,G,\widetilde{H},\widetilde{G}$, что $V(H)=\{x_1,...,x_k\}$,
$V(G)=\{x_1,...,x_l\}$,
$V(\widetilde{H})=\{\widetilde{x}_1,...,\widetilde{x}_k\}$,
$V(\widetilde{G})=\{\widetilde{x}_1,...,\widetilde{x}_l\}$, причем
$H\subset G$, $\widetilde{H}\subset\widetilde{G}$ (тем самым,
$k<l$). Граф $\widetilde{G}$ называется $(G,H)$--\emph{расширением
графа $\widetilde{H}$}, когда
$$
 \{x_{i_1},x_{i_2}\}\in E(G)\setminus E(H) \Rightarrow
 \{\widetilde{x}_{i_1},\widetilde{x}_{i_2}\}\in
 E(\widetilde{G})\setminus E(\widetilde{H}).
$$
Если выполняется соотношение
$$
 \{x_{i_1},x_{i_1}\}\in
 E(G)\setminus E(H) \Leftrightarrow
 \{\widetilde{x}_{i_1},\widetilde{x}_{i_2}\}\in
 E(\widetilde{G})\setminus E(\widetilde{H}),
$$
то $\widetilde{G}$ назовем \emph{точным расширением}, а пары
$(G,H)$ и $(\widetilde G,\widetilde H)$ --- изоморфными.
Зафиксируем число $\alpha>0$. Положим
$$
 v(G,H)=|V(G)\setminus V(H)|, \,\,
 e(G,H)=|E(G)\setminus E(H)|,
$$
$$
 f_{\alpha}(G,H)=v(G,H)-\alpha e(G,H).
$$
Если для любого такого графа $S,$ что $H\subset S\subseteq G,$
выполнено неравенство $f_{\alpha}(S,H)>0,$ то пара $(G,H)$
называется \emph{$\alpha$--надежной} (см. \cite{Janson},
\cite{Alon}, \cite{Spencer}). Если же для любого такого $S,$ что
$H\subseteq S\subset G,$ выполнено неравенство
$f_{\alpha}(G,S)<0,$ то пара $(G,H)$ называется {\it
$\alpha$--жесткой} (см. \cite{Janson}, \cite{Alon}). Если любая
вершина графа $H$ соединена с некоторой вершиной из $V(G)\setminus
V(H)$ и для любого такого графа $S$, что $H\subset S\subset G,$
справедливо неравенство $f_{\alpha}(S,H)>0$, но
$f_{\alpha}(G,H)=0$, то пару $(G,H)$ назовем
\emph{$\alpha$-нейтральной}. Введем, наконец, понятие максимальной
пары. Пусть $\widetilde{H}\subset\widetilde{G}\subset\Gamma$ и
$T\subset K$, причем $|V(T)|\leq|V(\widetilde{G})|.$ Пару
$(\widetilde{G},\widetilde{H})$ назовем
\emph{$(K,T)$--максимальной в $\Gamma$}, если у любого такого
подграфа $\widetilde{T}$ графа $\widetilde{G},$ что
$|V(\widetilde{T})|=|V(T)|$ и
$\widetilde{T}\cap\widetilde{H}\neq\widetilde{T},$ не существует
такого точного $(K,T)$-расширения $\widetilde{K}$ в
$\Gamma\setminus(\widetilde{G}\setminus\widetilde{T})$, что каждая
вершина из $V(\widetilde{K})\setminus V(\widetilde{T})$ не
соединена ребром ни с одной вершиной из $V(\widetilde{G})\setminus
V(\widetilde{T})$. Граф $\widetilde{G}$ назовем
\emph{$(K,T)$--максимальным в $\Gamma$}, если у любого такого
подграфа $\widetilde{T}$ графа $\widetilde{G},$ что
$|V(\widetilde{T})|=|V(T)|$, не существует такого точного
$(K,T)$-расширения $\widetilde{K}$ в
$\Gamma\setminus(\widetilde{G}\setminus\widetilde{T})$, что каждая
вершина из $V(\widetilde{K})\setminus V(\widetilde{T})$ не
соединена ребром ни с одной вершиной из $V(\widetilde{G})\setminus
V(\widetilde{T})$.

Обратимся теперь к случайному графу $G(N,p)$. Пусть
$\alpha\in(0,1]$, $p=N^{-\alpha}$. Пусть, кроме того, пара $(G,H)$
является $\alpha$--надежной и $V(H)=\{x_1,...,x_k\},$
$V(G)=\{x_1,...,x_l\}$. Рассмотрим произвольные вершины
$\widetilde{x}_1,...,\widetilde{x}_k\in V_N$ и случайную величину
$N_{(G,H)}(\widetilde{x}_1,...,\widetilde{x}_k)$ на вероятностном
пространстве $(\Omega_N,\mathcal{F}_N,{\sf P}_{N,p})$, которая
каждому графу $\mathcal{G}$ из $\Omega_N$ ставит в соответствие
количество $(G,H)$-расширений подграфа в $\mathcal{G}$,
индуцированного на $\{\widetilde{x}_1,...,\widetilde{x}_k\}$ (граф
$X$ является \emph{подграфом графа $Y$, индуцированным на
множество $S\subset V(Y)$}, если $V(X)=S$ и для любых вершин
$x,y\in S$ справедливо $\{x,y\}\in E(X)\Leftrightarrow\{x,y\}\in
E(Y)$). Иными словами, пусть $W\subset V_N
\setminus\{\widetilde{x}_1,...,\widetilde{x}_k\}$
--- множество, мощность $|W|$ которого равна $l-k.$ Если можно так
занумеровать элементы множества $W$ числами $k+1,k+2,...,l,$ что
граф $\mathcal{G}|_{\{\widetilde{x}_1,...,\widetilde{x}_l\}}$
является $(G,H)$--расширением графа
$\mathcal{G}|_{\{\widetilde{x}_1,...,\widetilde{x}_k\}},$ то
положим $I_W(\mathcal{G})=1.$ В противном случае
$I_W(\mathcal{G})=0.$ Случайная величина
$N_{(G,H)}(\widetilde{x}_1,...,\widetilde{x}_k)$ определяется
следующим равенством:
$$
N_{(G,H)}(\widetilde{x}_1,...,\widetilde{x}_k)=
\sum\limits_{W\subset V_N
\setminus\{\widetilde{x}_1,...,\widetilde{x}_k\},\, |W|=l-k}I_W.
$$

\begin{theorem} [\cite{Spencer}]
С асимптотической вероятностью $1$ для любых вершин
$\widetilde{x}_1,...,\widetilde{x}_k$ справедливо соотношение
$N_{(G,H)}(\widetilde{x}_1,...,\widetilde{x}_k)\sim {\sf E}_{N,p}
N_{(G,H)}(\widetilde{x}_1,...,\widetilde{x}_k).$ При этом
$$
 {\sf
 E}_{N,p}N_{(G,H)}(\widetilde{x}_1,...,\widetilde{x}_k))=\Theta(N^{f_{\alpha}(G,H)}).
$$
\label{safe_extensions}
\end{theorem}

Помимо этой теоремы Дж.~Спенсер и С.~Шела (см. \cite{Alon},
\cite{Shelah}) для исследования законов нуля или единицы доказали
теорему о количестве максимальных расширений подграфов в случайном
графе (для случая "запрещенных"$\,$жестких пар). В 2010 году мы
расширили их результат, рассмотрев "запрещенные"$\,$нейтральные
пары (см. \cite{zhuk_extensions}).

Более формально, пусть случайная величина
$N^{(K,T)}_{(G,H)}(\widetilde{x}_1,...,\widetilde{x}_k)$ ставит в
соответствие каждому графу $\mathcal{G}$ из $\Omega_N$ количество
таких точных $(G,H)$-расширений $\widetilde{G}$ графа
$\widetilde{H}=\mathcal{G}|_{\{\widetilde{x}_1,...,\widetilde{x}_k\}}$,
что пара $(\widetilde{G},\widetilde{H})$ является
$(K,T)$--максимальной в $\mathcal{G}$. Сформулируем теорему,
доказанную в работе Дж.~Спенсера и С.~Шела (см. \cite{Shelah}),
об асимптотическом поведении этой случайной величины.\\

\begin{theorem} [\cite{Shelah}]
Пусть пара $(K,T)$ является $\alpha$--жесткой. Тогда с
асимптотической вероятностью $1$ для любых вершин
$\widetilde{x}_1,...,\widetilde{x}_k$ выполнено
$$
 N^{(K,T)}_{(G,H)}(\widetilde{x}_1,...,\widetilde{x}_k)\sim
 N_{(G,H)}(\widetilde{x}_1,...,\widetilde{x}_k) \sim
 {\sf
 E}_{N,p}
 N^{(K,T)}_{(G,H)}(\widetilde{x}_1,...,\widetilde{x}_k)
 =\Theta\left(N^{f_{\alpha}(G,H)}\right).
$$
\label{maximal_extensions_rigid}
\end{theorem}



\subsection{Специфические графы}
\label{cicles}

Пусть $m\in\mathbb{N}$ --- произвольное натуральное число.
Рассмотрим такую пару графов $(G,H)$, что $G\supset H$. Будем
говорить, что граф $G$ является {\it $m$--расширением графа $H$
первого типа}, если $m\geq 3$ и выполнено следующее условие.
Существует такая вершина $x_1$ графа $G$, что
$$
 V(G)\setminus V(H)=\{y_1^1,...,y_{t_1}^1,y_1^2,...,y_{t_2}^2\},
$$
$$
 E(G)\setminus E(H)=\{\{x_1,y_1^1\},\{y_1^1,y_2^1\},...,\{y_{t_1-1}^1,y_{t_1}^1\},
 \{y_{t_1}^1,y_{1}^2\},\{y_1^2,y_2^2\},...,\{y_{t_2-1}^2,y_{t_2}^2\},\{y_{t_2}^2,y_{t_1}^1\}\},
$$
где $t_1+t_2\leq m-1$, $t_1\geq 0$, $t_2\geq 2$  и
$\rho^{\max}(G)<\frac{m}{m-1}$ (при $t_1=0$ вершина $x_1$ является
смежной с вершинами $y^2_1,y^2_{t_2}$). Граф $G$ мы называем {\it
$m$--расширением графа $H$ второго типа}, если $m\geq 2$ и
выполнено следующее условие. Существуют две такие различные
вершины $x_1,x_2$ графа $G$, что
$$
 G=(V(H)\sqcup\{y_1,...,y_t\},E(H)\sqcup\{\{x_1,y_1\},\{y_1,y_2\},...,\{y_{t-1},y_t\},\{y_t,x_2\}\}),
$$
где $t\leq m-1$ и $\rho^{\max}(G)<\frac{m}{m-1}$. Граф $G$
является {\it $m$--расширением графа $H$ третьего типа}, если
$m\geq 2$, $V(H)=V(G)$, $E(H)\subset E(G)$ и
$\rho^{\max}(G)<\frac{m}{m-1}$.

Для произвольного натурального числа $m\geq 3$ определим множество
графов $\mathcal{H}_m$. Пусть $x$ --- вершина. Граф без ребер на
множестве вершин $\{x\}$, принадлежит $\mathcal{H}_m$. Далее пусть
$G\in\mathcal{H}_m$. Множество $\mathcal{H}_m$ содержит все
попарно неизоморфные $m$--расширения первого, второго и третьего
типов графа $G$.

Заметим, что любой граф $G$ из $\mathcal{H}_m$, отличный от
$(\{x\},\varnothing)$, содержит в себе конечный набор таких
вложенных графов $G_1,...,G_t$, $G_0=(\{x\},\varnothing)\subset
G_1\subset...\subset G_t\subseteq G$, что выполнены следующие
свойства:

\begin{itemize}

\item[---] $G_i\neq G_{i+1}$ для всех $i\in\{0,1,...,t-1\}$,

\item[---] граф $G$ либо совпадает с $G_t$, либо является
$m$--расширением третьего типа графа $G_t$, а графы $G_i$ являются
$m$--расширениями первого или второго типа графов $G_{i-1}$ при
$i\in\{0,1,...,t\}$.

\end{itemize}

Такую последовательность графов $G_0,G_1,...,G_t,G$ будем называть
{\it $m$--разложением} графа $G$.

Сформулируем утверждение о свойствах множества $\mathcal{H}_m$.

\begin{lemma} Выполнены следующие свойства.

\begin{itemize}

\item[1.] Пусть $G\in\mathcal{H}_m$ и $G_0,G_1,...,G_t,G$ --- его
$m$--разложение, причем либо $t=1$ и $G_t\neq G$, либо $t\geq 2$.
Тогда найдутся такие натуральные числа $a,b$, что $b\leq m$ и
$\rho^{\max}(G)=1+\frac{1}{m-1+b/a}$.


\item[2.] Пусть $m\geq 2$ и $\rho\in(1,m/(m-1))$ --- произвольное
число. Тогда существует такое число $\eta\in\mathbb{N}$, что для
любого натурального $v>\eta$ и у любого графа $G\in\mathcal{H}_m$
на $v$ вершинах найдется подграф на не более чем $\eta$ вершинах с
плотностью, превосходящей $\rho$.

\end{itemize}

\label{H-properties}
\end{lemma}

{\bf Доказательство} начнем со свойства 1. Если $t=1$, то
$\rho(G_1)=1,$ $v(G_1)\leq m$. Следовательно,
$\rho(G)=\frac{v(G_1)+e}{v(G_1)}$, где $e\in\mathbb{N}$. Поэтому
$\rho(G)=1+\frac{1}{1+\frac{v(G_1)-e}{e}}$. Так как $v(G_1)-e<m$,
то свойство 1 для рассмотренного случая доказано.

Пусть $t\geq 2$. Положим $v^{i+1}=v(G_{i+1})-v(G_i)$,
$i\in\{0,1,...,t-1\}$. Докажем, что
$\rho(G_t)=1+\frac{1}{m-1+b_1/(t-1)}$ при некотором $b_1\leq m$.
Справедливы соотношения
$$
 \rho(G_t)=\frac{v^1+...+v^{t}+t}{1+v^1+...+v^{t}}=1+\frac{t-1}{1+v^1+...+v^{t}}<1+\frac{1}{m-1}.
$$

Так как $v^i\leq m-1$, $i\in\{1,...,t\}$, то
$$
 m-1<\frac{1+v^1+...+v^{t}}{t-1}\leq\frac{1+t(m-1)}{t-1}=m-1+\frac{m}{t-1}.
$$
Таким образом,
$$
 1+\frac{1}{m-1+\frac{m}{t-1}}\leq \rho(G)<1+\frac{1}{m-1},
$$
при этом
$$
 \rho(G)=1+\frac{1}{m-1+\frac{b_1}{t-1}},
$$
где $b_1=1+v^1+...+v^t-(m-1)(t-1)$. Поэтому число $b_1$ не
превосходит $m$.

Докажем теперь, что $\rho(G)=1+\frac{1}{m-1+b_2 / (t+e_0-1)},$ где
$b_2\leq m, e_0=e(G)-e(G_t)$. Имеем
$$
 \rho(G)=\frac{m(t-1)+b_1+e_0}{(m-1)(t-1)+b_1}=
 1+\frac{1}{m-1+\frac{b_1-(m-1)e_0}{t-1+e_0}}.
$$
Так как $\rho(G)<1+\frac{1}{m-1},$ то $0<b_2\leq b_1\leq m,$ где
$b_2=b_1-(m-1)e_0$.


Пусть, наконец $H\subset G,$ $\rho(G)<\rho(H)<\frac{m}{m-1}$.
Тогда
$$
 \rho(H)=\frac{m(t-1)+b_1+e_0-y}{(m-1)(t-1)+b_1-x}
$$
для некоторых натуральных чисел $x,y$. Справедливо неравенство
$y\geq x,$ так как граф $G$
--- связанный. Докажем, что $\rho(H)=1+\frac{1}{m-1+b /
(t+e_0-1+x-y)},$ где $b\leq m$. Имеем
$$
 \rho(H)=1+\frac{1}{m-1+\frac{b_1+y(m-1)-mx-(m-1)e_0}{t-1+e_0+x-y}}.
$$

Так как $\rho(H)>\rho(G),$ то
$$
 \frac{b_1+y(m-1)-mx-(m-1)e_0}{t-1+e_0+x-y}<\frac{b_1-(m-1)e_0}{t-1+e_0}.
$$
Но знаменатель первой дроби не больше, чем знаменатель второй.
Следовательно, $b=b_1+y(m-1)-mx-(m-1)e_0<b_1-(m-1)e_0\leq m,$ что
и требовалось доказать.\\

Перейдем к доказательству свойства 2. В силу определения
$m$--расширений первого и второго типа, если $G_0,G_1,...,G_t,G$
--- $m$--разложение некоторого графа $G\in\mathcal{H}_m$, то
$v(G_{i+1})-v(G_i)\leq m-1,$ $e(G_{i+1})-e(G_i)=v+1$ для всех
$i\in\{0,...,t-1\}$. Следовательно, для каждого $n\in\mathbb{N}$ и
для любого графа $G\in\mathcal{H}_m$, количество вершин которого
не меньше, чем $(m-1)n+1$, справедливо неравенство
$\rho(G)\geq\frac{mn}{(m-1)n+1}$. Действительно, плотность любого
графа из множества $\mathcal{H}_m$ меньше, чем $\frac{m}{m-1}$,
поэтому при добавлении к нему его $m$--расширения первого или
второго типа его плотность увеличивается. Кроме того, для любого
$\rho\in(1,\frac{m}{m-1})$ найдется такое число
$n_0\in\mathbb{N}$, что при всех натуральных $n\geq n_0$ выполнено
неравенство $\frac{mn}{(m-1)n+1}>\rho.$ Заметим, наконец, что в
любом графе из $\mathcal{H}_m$ на более чем $(m-1)(n+1)+1$
вершинах найдется подграф из $\mathcal{H}_m$, количество вершин
которого находится в отрезке $[(m-1)n+1,(m-1)(n+1)+1]$. Поэтому,
очевидно, для $\eta=(m-1)(n_0+1)+1$ утверждение леммы выполнено.

\section{Игра Эренфойхта}
\label{game}

Основным средством в доказательстве законов нуля или единицы для
свойств первого порядка случайных графов, как мы заметили выше,
служит теорема А. Эренфойхта, доказанная в 1960 году (см.
\cite{Ehren}). В данном разделе мы сформулируем ее частный случай
для графов. Прежде всего определим игру Эренфойхта $\EHR(G,H,i)$
на двух графах $G,H$ с количеством раундов, равным $i$ (см.,
например, \cite{Janson}, \cite{Alon}). Пусть
$V(G)=\{x_1,...,x_n\},$ $V(H)=\{y_1,...,y_m\}$. В $\nu\mbox{-}$ом
раунде ($1 \leq \nu \leq i$) Новатор выбирает вершину из любого
графа (он выбирает либо $x_{j_{\nu}}\in V(G)$, либо
$y_{j'_{\nu}}\in V(H)$). Затем Консерватор выбирает вершину из
оставшегося графа. Если Новатор выбирает в $\mu\mbox{-}$ом раунде,
скажем, вершину $x_{j_{\mu}}\in V(G),$ $j_{\mu}=j_{\nu}$
($\nu<\mu$), то Консерватор должен выбрать $y_{j'_{\nu}}\in V(H)$.
Если же в этом раунде Новатор выбирает, скажем, вершину
$x_{j_{\mu}}\in V(G),$ $j_{\mu}\notin\{j_1,...,j_{\mu-1}\},$ то и
Консерватор должен выбрать такую вершину $y_{j'_{\mu}}\in V(H)$,
что $j'_{\mu}\notin\{j'_1,...,j'_{\mu-1}\}.$ Если он не может
этого сделать, то игру выигрывает Новатор. К концу игры выбраны
вершины $x_{j_1},...,x_{j_{i}}\in V(G)$, а также вершины
$y_{j'_1},...,y_{j'_{i}}\in V(H)$. Некоторые из этих вершин могут
совпадать. Выберем из них только различные: $x_{h_1},...,x_{h_l};$
$y_{h'_1},...,y_{h'_l},$ $l \leq i.$ Консерватор побеждает тогда и
только тогда, когда соответствующие подграфы изоморфны с точностью
до порядка вершин:
$$
 G|_{\{x_{h_1},...,x_{h_l}\}}\cong
 H|_{\{y_{h'_1},...,y_{h'_l}\}}.
$$

\begin{theorem} [\cite{Ehren}]
Для любых двух графов $G,H$ и любого $i\in\mathbb{N}$ Консерватор
имеет выигрышную стратегию в игре $\EHR(G,H,i)$ тогда и только
тогда, когда для любого свойства $L$ первого порядка, выражаемого
формулой, кванторная глубина которой не превышает $i$, либо оба
графа обладают этим свойством, либо оба не обладают.
\end{theorem}

Несложно показать, что из этой теоремы вытекает следующее
следствие о законах нуля или единицы (см., например,
\cite{Janson}, \cite{Alon}). Пусть ${\sf P}_{N,p(N)}\times {\sf
P}_{M,p(M)}$
--- прямое произведение мер ${\sf P}_{N,p(N)}, {\sf P}_{M,p(M)}$.

\begin{theorem}
Случайный граф $G(N,p)$ подчиняется $k$-закону нуля или единицы
тогда и только тогда
$$
 \lim\limits_{N,M\rightarrow\infty}{\sf P}_{N,p(N)}\times {\sf
 P}_{M,p(M)}(\mbox{у Консерватора есть выигрышная стратегия в
 игре}
$$
$$
 \EHR(G(N,p(N)),G(M,p(M)),k))=1.
$$

\label{ehren}
\end{theorem}

Обратимся к описанию выигрышной стратегии Консерватора. Всюду
далее мы считаем, что в каждом раунде игроки выбирают вершины
отличные от уже выбранных. Такое предположение не ограничивает
общности в силу того, что размеры графов в рассматриваемых случаях
а также утверждения, используемые в рассуждениях, позволяют
выбирать вершины, отличные от уже выбранных. Сначала мы опишем
некоторые свойства графов $G,H$, благодаря которым Консерватор
выигрывает, придерживаясь такой стратегии. В следующем разделе мы
докажем, что случайный граф обладает этими свойствами с
вероятностью, стремящейся к $1$.

Пусть $n_1,n_2,n_3,n_4$ --- некоторые натуральные числа, $n_2\leq
n_1$, $n_4\leq n_3$, $\rho$
--- произвольное положительное число. Будем говорить, что граф $G$
является {\it $(n_1,n_2,n_3,n_4,\rho)$--разреженным}, если он
обладает следующими свойствами.

\begin{itemize}

\item[1)] Пусть $K$ --- граф, количество вершин которого не
превосходит $n_1$. Если $\rho^{\max}(K)<\rho$, то $G$ содержит
подграф, изоморфный $K$. Если $\rho^{\max}(K)>\rho$, то $G$ не
содержит подграф, изоморфный $K$.

\item[2)] Пусть $\mathcal{H}$ --- множество таких
$1/\rho$--надежных пар $(H_1,H_2)$, что $v(H_1)\leq n_1$,
$v(H_2)\leq n_2.$ Пусть $\mathcal{K}$ --- множество таких
$1/\rho$--жестких пар $(K_1,K_2)$, что $v(K_1)\leq n_3$,
$v(K_2)\leq n_4.$ Тогда для любых пар $(H_1,H_2)\in\mathcal{H},$
$(K_1,K_2)\in\mathcal{K}$ и для любого подграфа $G_2\subset G$ на
$v(H_2)$ вершинах в графе $G$ найдется подграф $G_1$, являющийся
$(K_1,K_2)$--максимальным в $G$ точным $(H_1,H_2)$--расширением
графа $G_2$.

\end{itemize}

Выигрышная стратегия Консерватора, описанная ниже, опирается
именно на свойство $(n_1,n_2,n_3,n_4,\rho)$--разреженности обоих
графов, на которых играют Новатор и Консерватор (при некоторых
значениях $n_1,n_2,n_3,n_4,\rho$).

Пусть $\rho\in(1,\frac{2^{k-1}}{2^{k-1}-1})$,
$\rho\notin\{1+\frac{1}{2^{k-1}-1+b/a}, \, a,b\in\mathbb{N}, b\leq
2^{k-1}\}$. При доказательстве закона нуля или единицы в следующем
разделе выбор числа $\rho$ будет зависеть от числа $\alpha$, для
которого доказывается закон ($\rho=1/\alpha$). Такое число $\rho$
является порогом для значений плотностей подграфов и расширений
(их можно определить по аналогии с плотностями графов),
содержащихся в случайном графе. Описанная ниже стратегия опирается
на выбор числа $\rho$. Таким образом, при доказательстве закона
нуля или единицы для случайного графа $G(N,N^{-\alpha})$ мы
сначала фиксируем число
$\rho(\alpha)\in(1,\frac{2^{k-1}}{2^{k-1}-1})$, затем выбираем
числа $n_1,n_2,n_3,n_4,$ некоторым образом зависящие от $\rho$ и
$k$, доказываем, что случайный граф с вероятностью, стремящейся к
1, является $(n_1,n_2,n_3,n_4,\rho)$--разреженным, после чего
пользуемся стратегией, описанной ниже.

Обозначим $\eta(\rho)$ число из формулировки леммы 1, т.е. такое
число, что для любого натурального $v>\eta(\rho)$ и у любого графа
$G\in\mathcal{H}_{2^{k-1}}$ на $v$ вершинах найдется подграф на не
более чем $\eta(\rho)$ вершинах с плотностью, превосходящей
$\rho$. Положим
\begin{equation}
 n_1(\rho)=\eta(\rho)+(k-1)\left(\left[\frac{1}{\rho-1}\right]+1\right), \quad
 n_2(\rho)=\eta(\rho)+(k-2)\left(\left[\frac{1}{\rho-1}\right]+1\right),
\label{n1}
\end{equation}
\begin{equation}
 n_3=2^{k-2}+1, \quad n_4=2.
\label{n4}
\end{equation}
Пусть графы $G,H$ являются
$(n_1(\rho),n_2(\rho),n_3,n_4,\rho)$--разреженными. Опишем
стратегию Консерватора в игре $\EHR(G,H,k)$. Будем обозначать
$X_i$ граф, выбранный Новатором в $i$-ом раунде. Оставшийся граф
будем обозначать $Y_i$. Вершины, выбранные в графе $X_i$ в первых
$i$ раундах обозначим $x_i^1,...,x_i^i$, в графе $Y_i$ ---
$y_i^1,...,y_i^i$. Итак, пусть в первом раунде Новатор выбрал
вершину $x_1^1$. В силу леммы 1 и свойства
$(n_1(\rho),n_2(\rho),n_3,n_4,\rho)$--разреженности графа $X_1$ в
нем не существует подграфа, изоморфного некоторому графу из
$\mathcal{H}_{m}$ с количеством вершин, превосходящим
$\eta(\rho).$ Обозначим $\widetilde{X}_1^1$ подграф в $X_1$ на
$v_1$ вершинах, изоморфный некоторому графу из
$\mathcal{H}_{2^{k-1}}$, содержащий вершину $x_1^1$ и обладающий
следующим свойством. В $X_1$ не существует подграфа, содержащего
вершину $x_1^1$ и изоморфного некоторому графу из
$\mathcal{H}_{2^{k-1}},$ количество вершин которого превосходит
$v_1$. Выполнены неравенства $v_1\leq \eta(\rho)<n_1(\rho)$.
Поэтому в силу леммы 1 и свойства
$(n_1(\rho),n_2(\rho),n_3,n_4,\rho)$--разреженности графа $X_1$
плотность графа $\widetilde{X}_1$ меньше, чем $\rho$.
Следовательно, по свойству
$(n_1(\rho),n_2(\rho),n_3,n_4,\rho)$--разреженности графа $Y_1$ в
нем найдется подграф $\widetilde{Y}_1^1$, изоморфный
$\widetilde{X}_1$. Пусть при соответствующем изоморфизме
$\varphi_1:\widetilde{X}_1\rightarrow\widetilde{Y}_1^1$ вершина
$x_1^1$ переходит в вершину $y_1^1$, которую и выберет Консерватор
в первом раунде.

Пусть сыграно $i$ раундов, $1\leq i < k$. Опишем стратегию
Консерватора в $i+1$-ом раунде. Ниже мы определим ряд свойств
подграфов в $G,H$ (они обозначены (I), (II), (III)). В
соответствии со свойством (I) выбранные вершины
$x_1^i,...,x_i^i,y_1^i,...,y_i^i$ должны принадлежать объединению
этих подграфов. Мы предположим, что в $G$ и $H$ содержатся
подграфы, обладающие этими свойствами. Затем мы докажем, что
независимо от выбора Новатором вершины $x_{i+1}^{i+1}$ Консерватор
сможет найти такую вершину $y_{i+1}^{i+1}$, что вершины
$x_1^{i+1},...,x_{i+1}^{i+1},y_1^{i+1},...,y_{i+1}^{i+1}$ также
будут содержаться в подграфах, обладающих свойствами (I), (II) и
(III). Кроме того, мы докажем, что вершины $x_1^1,y_1^1$
содержатся в подграфах, обладающих упомянутыми свойствами, откуда
по индукции последует аналогичное утверждение для последнего
раунда, т.е. для вершин $x_1^k,...,x_k^k,y_1^k,...,y_k^k$. В
частности, из этих свойств мы выведем, что графы
$X_k|_{\{x_1^k,...,x_k^k\}}$ и $Y_k|_{\{y_1^k,...,y_k^k\}}$
изоморфны.

Введем обозначение, которое мы будем использовать в дальнейшем в
этом разделе. Для произвольного графа $Q$, его подграфа $W$ и
вершины $x\in V(Q)\setminus V(W)$ обозначим $d_Q(x,W)$ длину самой
короткой цепи из $Q$ (количество ребер в ней), соединяющую вершину
$x$ с произвольной вершиной графа $W$ (очевидно, все вершины такой
цепи кроме одной не принадлежат графу $W$). Для подграфов
$W_1,W_2\subset Q$, не имеющих общих вершин, положим
$d_Q(W_1,W_2)=\min_{x\in V(W_1)}d_Q(x,W_2)$. Пусть $r$
--- произвольное натуральное число, не превосходящее $i$. Пусть,
кроме того, $W_1,...,W_r$ --- подграфы в $Q$. Будем говорить, что
$W_1,...,W_r$ обладают $(k,i,r)$--свойством в $Q$, если

\begin{itemize}

\item[---] любые два графа из $W_1,...,W_r$ не имеют общих вершин;

\item[---] для любых различных $j_1,j_2\in\{1,...,r\}$ справедливо
неравенство $d_{Q}(W_{j_1},W_{j_2})>2^{k-i}$;


\item[---] для любого $j\in\{1,...,r\}$ в графе $Q$ не существует
подграфа, являющегося $2^{k-i}$--расширением первого или второго
типа графа $W_j$;


\item[---] мощность множества $|V(W_1\cup...\cup W_r)|$ не
превосходит
$\eta(\rho)+(i-1)\left(\left[\frac{1}{\rho-1}\right]+1\right)$.

\end{itemize}

Предположим, что для некоторого $r\in\{1,...,i\}$ граф $X_i$
содержит подграфы $\widetilde{X}_i^1,...,\widetilde{X}_i^r$, граф
$Y_i$ содержит подграфы $\widetilde{Y}_i^1,...,\widetilde{Y}_i^r$,
которые обладают следующими свойствами.

\begin{itemize}

\item[(I)] Вершины $x_i^1,...,x_i^i$ принадлежат множеству
$V(\widetilde{X}_i^1\cup...\cup \widetilde{X}_i^r)$, вершины
$y_i^1,...,y_i^i$ принадлежат множеству
$V(\widetilde{Y}_i^1\cup...\cup \widetilde{Y}_i^r)$.

\item[(II)] Графы $\widetilde{X}^1_i,...,\widetilde{X}^r_i$
обладают $(k,i,r)$--свойством в $X_i$, графы
$\widetilde{Y}^1_i,...,\widetilde{Y}^r_i$ обладают
$(k,i,r)$--свойством в $Y_i$.

\item[(III)] Графы $\widetilde{X}^j_i$ и $\widetilde{Y}^j_i$
изоморфны при каждом $j\in\{1,...,r\}$ и при некотором
соответствующем изоморфизме (общим для всех графов, так как они не
имеют общих вершин) вершины $x_i^j$ переходят в вершины $y_i^j,$
$j\in\{1,...,i\}$.

\end{itemize}

Если $X_{i+1}=X_i$, то положим
$\widetilde{X}^j_{i+1}=\widetilde{X}^j_i$,
$\widetilde{Y}^j_{i+1}=\widetilde{Y}^j_i$, $j\in\{1,...,r\}$. В
противном случае положим
$\widetilde{X}^j_{i+1}=\widetilde{Y}^j_i$,
$\widetilde{Y}^j_{i+1}=\widetilde{X}^j_i$, $j\in\{1,...,r\}$.
Пусть $\varphi_{i+1}$ --- изоморфизм из
$\widetilde{X}^1_{i+1}\cup...\cup\widetilde{X}^r_{i+1}$ в
$\widetilde{Y}^1_{i+1}\cup...\cup\widetilde{Y}^r_{i+1}$,
переводящий графы $\widetilde{X}^j_{i+1}$ в
$\widetilde{Y}^j_{i+1}$ при $j\in\{1,...,r\}$. Пусть, кроме того,
$\varphi_{i+1}(x^j_{i+1})=y^j_{i+1}$ при $j\in\{1,...,i\}$.
Рассмотрим далее три различные ситуации.

\begin{enumerate}

\item Предположим, что Новатор в $i+1$-ом раунде выбрал вершину
$x_{i+1}^{i+1}$ из множества
$V(\widetilde{X}^1_{i+1}\cup...\cup\widetilde{X}^r_{i+1})$. Тогда
Консерватор выберет $y_{i+1}^{i+1}=\varphi(x_{i+1}^{i+1})$.
Заметим, что мы определили графы
$\widetilde{X}_{i+1}^1,...,\widetilde{X}_{i+1}^r,\widetilde{Y}_{i+1}^1,...,\widetilde{Y}_{i+1}^r$
просто переобозначив графы
$\widetilde{X}_{i}^1,...,\widetilde{X}_{i}^r,\widetilde{Y}_{i}^1,...,\widetilde{Y}_{i}^r$
и не меняя их структуры. Поэтому, как нетрудно видеть, при $i<k-1$
графы $\widetilde{X}^1_{i+1},...,\widetilde{X}^r_{i+1}$ обладают
$(k,i+1,r)$--свойством в $X_{i+1}$, графы
$\widetilde{Y}^1_{i+1},...,\widetilde{Y}^r_{i+1}$ обладают
$(k,i+1,r)$--свойством в $Y_{i+1}$. Кроме того, графы
$\widetilde{X}^j_{i+1}$ и $\widetilde{Y}^j_{i+1}$ изоморфны при
каждом $j\in\{1,...,r\}$ и при соответствующем изоморфизме
$\varphi_{i+1}$ (одном и том же для всех графов) вершины
$x_{i+1}^j$ переходят в вершины $y_{i+1}^j,$ $j\in\{1,...,i+1\}$.
Иными словами, для $i+1$-го раунда мы подобрали графы
$\widetilde{X}^1_{i+1},...,\widetilde{X}^r_{i+1},\widetilde{Y}^1_{i+1},...,\widetilde{Y}^r_{i+1}$,
обладающие свойствами (I), (II), (III).

\item Предположим теперь, что Новатор выбрал вершину
$x_{i+1}^{i+1}$, не принадлежащую множеству
$V(\widetilde{X}^1_{i+1}\cup...\cup\widetilde{X}^r_{i+1})$, при
этом
$d_{X_{i+1}}(x_{i+1}^{i+1},\widetilde{X}_{i+1}^1\cup...\cup\widetilde{X}_{i+1}^r)\leq
2^{k-1-i}$. Заметим, что в силу определения графов
$\widetilde{X}_{i+1}^1,...,\widetilde{X}_{i+1}^r$ в графе
$X_{i+1}$ найдется ровно одна цепь
$c_{X_{i+1}}$,
проходящая только через вершины графа
$X_{i+1}\setminus\left(\widetilde{X}_{i+1}^1\cup...\cup\widetilde{X}_{i+1}^r\right)$
(не считая последней вершины), длины, не превосходящей
$2^{k-1-i}$, соединяющая $x_{i+1}^{i+1}$ с некоторой вершиной
$\widetilde{x}_{i+1}^l$ графа
$\widetilde{X}_{i+1}^1\cup...\cup\widetilde{X}_{i+1}^r$, где
$l\in\{1,...,r\}$, $\widetilde{x}_{i+1}^l\in
V(\widetilde{X}_{i+1}^l)$. Следовательно, пара
$(\widetilde{X}_{i+1}^1\cup...\cup\widetilde{X}_{i+1}^r\cup
 c_{X_{i+1}},
 \widetilde{X}_{i+1}^1\cup...\cup\widetilde{X}_{i+1}^r)$ является $1/\rho$--надежной. Кроме того,
$$
 |V(\widetilde{Y}_{i+1}^1\cup...\cup\widetilde{Y}_{i+1}^r)|<\eta(\rho)+(i-1)\left(\left[\frac{1}{\rho-1}\right]+1\right)\leq
 \eta(\rho)+(k-2)\left(\left[\frac{1}{\rho-1}\right]+1\right).
$$
Поэтому в силу свойства
$(n_1(\rho),n_2(\rho),n_3,n_4,\rho)$--разреженности графа
$Y_{i+1}$ в нем найдется точное $(K_1,K_2)$--максимальное
$(\widetilde{X}_{i+1}^1\cup...\cup\widetilde{X}_{i+1}^r\cup
c_{X_{i+1}},
\widetilde{X}_{i+1}^1\cup...\cup\widetilde{X}_{i+1}^r)$--расширение
графа $\widetilde{Y}_{i+1}^1\cup...\cup\widetilde{Y}_{i+1}^r$ для
всех таких $1/\rho$--жестких пар $(K_1,K_2)$, что $v(K_2)=2$,
$v(K_1)\leq 2^{k-1-i}$. Действительно,
$$
 |V(\widetilde{X}_{i+1}^1\cup...\cup\widetilde{X}_{i+1}^r\cup
 c_{X_{i+1}})|\leq
 \eta(\rho)+(k-2)\left(\left[\frac{1}{\rho-1}\right]+1\right)+2^{k-1-i}<
$$
$$
 <\eta(\rho)+(k-2)\left(\left[\frac{1}{\rho-1}\right]+1\right)+\frac{1}{\rho-1}\leq
 \eta(\rho)+(k-1)\left(\left[\frac{1}{\rho-1}\right]+1\right)=n_1(\rho).
$$
Иными словами, существует такая вершина $y_{i+1}^{i+1}\in
V(Y_{i+1})$, что
$$
 d_{Y_{i+1}}(y_{i+1}^{i+1},\widetilde{Y}_{i+1}^l)=
 d_{X_{i+1}}(x_{i+1}^{i+1},\widetilde{X}_{i+1}^l),
$$
и единственная цепь $c_{Y_{i+1}}$ 
длины, не превосходящей $2^{k-1-i}$, соединяющая $y_{i+1}^{i+1}$ с
некоторой вершиной $\widetilde{y}_{i+1}^l$ графа
$\widetilde{Y}_{i+1}^1\cup...\cup\widetilde{Y}_{i+1}^r$,
$\widetilde{y}_{i+1}^l\in V(\widetilde{Y}_{i+1}^l)$. Переопределим
графы $\widetilde{X}_{i+1}^l,\widetilde{Y}_{i+1}^l:$
$$
 \widetilde{X}_{i+1}^l:=\widetilde{X}_{i+1}^l\cup
 c_{X_{i+1}},\,\,\,\,
 \widetilde{Y}_{i+1}^l:=\widetilde{Y}_{i+1}^l\cup c_{Y_{i+1}}.
$$
Остальные графы
$\widetilde{X}_{i+1}^1,...,\widetilde{X}_{i+1}^{l-1},\widetilde{X}_{i+1}^{l+1},...,\widetilde{X}_{i+1}^{r}$,
$\widetilde{Y}_{i+1}^1,...,\widetilde{Y}_{i+1}^{l-1},\widetilde{Y}_{i+1}^{l+1},...,\widetilde{Y}_{i+1}^{r}$
оставим без изменений. Продолжим изоморфизм графов $\varphi_{i+1}$
на вершины из множества $V(\widetilde{X}_{i+1}^l)$: для вершины
$v$ цепи
$c_{X_{i+1}}$ 
с номером $h$, считая от $\widetilde{x}_{i+1}^l$, найдем вершину
$u$ цепи
$c_{Y_{i+1}}$
с тем же номером $h$, считая от $\widetilde{y}_{i+1}^l$, и
определим $\varphi_{i+1}(v)=u$. Тогда
$\varphi_{i+1}|_{\widetilde{X}^j_{i+1}}:\widetilde{X}^j_{i+1}\rightarrow\widetilde{Y}^j_{i+1}$
--- изоморфизм при каждом $j\in\{1,...,r\}$ и
$\varphi_{i+1}(x_{i+1}^j)=y_{i+1}^j$ при всех $j\in\{1,...,i+1\}$,
т.е. графы
$\widetilde{X}_{i+1}^1,...,\widetilde{X}_{i+1}^r,\widetilde{Y}_{i+1}^1,...,\widetilde{Y}_{i+1}^r$
обладают свойствами (I) и (III). Докажем, что при $i<k-1$
выполнено свойство (II) (графы
$\widetilde{X}^1_{i+1},...,\widetilde{X}^r_{i+1}$ обладают
$(k,i+1,r)$--свойством в $X_{i+1}$, а графы
$\widetilde{Y}^1_{i+1},...,\widetilde{Y}^r_{i+1}$ обладают
$(k,i+1,r)$--свойством в $Y_{i+1}$). Для этого достаточно доказать
следующие утверждения:

\begin{itemize}

\item[---] для любого $j\in\{1,...,r\}\setminus\{l\}$
$$
 d_{X_{i+1}}(\widetilde{X}^{j}_{i+1},\widetilde{X}^{l}_{i+1})>2^{k-i-1},\,\,\,\,\,
 d_{Y_{i+1}}(\widetilde{Y}^{j}_{i+1},\widetilde{Y}^{l}_{i+1})>2^{k-i-1};
$$

\item[---] в графе $X_{i+1}$ (в графе $Y_{i+1}$) не найдется
подграфа, являющегося $2^{k-i-1}$--расширением графа
$\widetilde{X}_{i+1}^l$ (графа $\widetilde{Y}_{i+1}^l$) первого
или второго типа;


\item[---] мощность множеств
$|V(\widetilde{X}^1_{i+1}\cup...\cup\widetilde{X}^r_{i+1})|,\,|V(\widetilde{Y}^1_{i+1}\cup...\cup\widetilde{Y}^r_{i+1})|$
не превосходит
$\eta(\rho)+i\left(\left[\frac{1}{\rho-1}\right]+1\right)$.

\end{itemize}

Предположим, что найдется $j\in\{1,...,r\}\setminus\{l\}$ и цепь
длины, не превосходящей $2^{k-1-i}$, соединяющая некоторую вершину
$u$ графа $\widetilde{X}^l_{i+1}$ с некоторой вершиной $v$ графа
$\widetilde{X}^j_{i+1}$. Так как $2^{k-1-i}<2^{k-i}$, то $u\in
V(c_{X_{i+1}})$. 
Но длина цепи
$c_{X_{i+1}}$
не превосходит $2^{k-1-i}$. Следовательно,
$d_{X_{i+1}}(\widetilde{X}^j_{i+1},\widetilde{X}^l_{i+1}\setminus
(c_{X_{i+1}}
\setminus\{\widetilde{x}_{i+1}^l\}))\leq
2^{k-1-i}+2^{k-1-i}=2^{k-i}$. Но величина в левой части
неравенства равна либо
$d_{X_i}(\widetilde{X}^l_i,\widetilde{X}^j_i)$, либо
$d_{Y_i}(\widetilde{Y}^l_i,\widetilde{Y}^j_i)$. Получили
противоречие либо с $(k,i,r)$--свойством графов
$\widetilde{X}_i^1,...,\widetilde{X}_i^r$, либо с
$(k,i,r)$--свойством графов
$\widetilde{Y}_i^1,...,\widetilde{Y}_i^r$. Таким образом,
$d_{X_{i+1}}(\widetilde{X}^{j}_{i+1},\widetilde{X}^{l}_{i+1})>2^{k-i-1}$.
Аналогично доказывается неравенство
$d_{Y_{i+1}}(\widetilde{Y}^{j}_{i+1},\widetilde{Y}^{l}_{i+1})>2^{k-i-1}$.

Доказательство второго утверждения мы тоже приводим только для
графа $X_{i+1}$, так как оно в точности повторяет доказательство
для графа $Y_{i+1}$. Итак, пусть в графе $X_{i+1}$ существует
подграф $W$, являющийся $2^{k-1-i}$--расширением первого или
второго типа графа $\widetilde{X}_{i+1}^l$. Рассмотрим множество
ребер $E=E(W)\setminus (E(\widetilde{X}_{i+1}^l)\cup E(W\setminus
\widetilde{X}_{i+1}^l))$. Вершин, принадлежащих множеству
$V(\widetilde{X}_{i+1}^l)$ и являющихся концами ребер из $E$, не
более двух. Обозначим их $v_1$ и $v_2$ (вообще говоря, эти вершины
могут совпадать). Если $v_1,v_2\in
V(\widetilde{X}_{i+1}^{l})\setminus
(V(c_{X_{i+1}})\setminus\{\widetilde{x}_{i+1}^l\})$, то мы
приходим к противоречию либо с $(k,i,r)$--свойством графов
$\widetilde{X}_i^1,...,\widetilde{X}_i^r$, либо с
$(k,i,r)$--свойством графов
$\widetilde{Y}_i^1,...,\widetilde{Y}_i^r$. Если же хотя бы одна из
вершин $v_1,v_2$ не принадлежит множеству
$V(\widetilde{X}_{i+1}^{l})\setminus
(V(c_{X_{i+1}})\setminus\{\widetilde{x}_{i+1}^l\})$, то в графе
$W$ найдется подграф $W_1$, множество вершин которого содержит
$V(\widetilde{X}_{i+1}^{l})\setminus
(V(c_{X_{i+1}})\setminus\{\widetilde{x}_{i+1}^l\})$, являющийся
$2^{k-i}$--расширением первого или второго типа графа
$W|_{V(\widetilde{X}_{i+1}^{l})\setminus
(V(c_{X_{i+1}})\setminus\{\widetilde{x}_{i+1}^l\})}$. Мы снова
приходим к противоречию либо с $(k,i,r)$--свойством графов
$\widetilde{X}_i^1,...,\widetilde{X}_i^r$, либо с
$(k,i,r)$--свойством графов
$\widetilde{Y}_i^1,...,\widetilde{Y}_i^r$. Последнее утверждение
выполнено, так как количество добавленных вершин не превосходит
$2^{k-1-i}\leq \left[\frac{1}{\rho-1}\right]+1$.

\end{enumerate}

Если $i=k-1$, то в обоих случаях так как $\varphi_k:
\widetilde{X}_k^1\cup...\cup\widetilde{X}_k^r\rightarrow
\widetilde{Y}_k^1\cup...\cup\widetilde{Y}_k^r$ --- изоморфизм и
$\varphi_k(x_k^j)=y_k^j$ при всех $j\in\{1,...,k\}$, то графы
$X_k|_{\{x_k^1,...,x_k^k\}}$, $Y_k|_{\{y_k^1,...,y_k^k\}}$ также
изоморфны и Консерватор побеждает.

\begin{enumerate}

\item[3.] Пусть, наконец,
$d_{X_{i+1}}(x_{i+1}^{i+1},\widetilde{X}_{i+1}^1\cup...\cup\widetilde{X}_{i+1}^r)>
2^{k-1-i}$. Найдем подграф в $X_{i+1}$, содержащий наибольшее
количество вершин, одна из которых совпадает с $x_{i+1}^{i+1}$, и
изоморфный некоторому графу из множества
$\mathcal{H}_{2^{k-1-i}}$. Обозначим полученный граф
$\widetilde{X}_{i+1}^{r+1}$. В силу свойства
$(n_1(\rho),n_2(\rho),n_3,n_4,\rho)$--разреженности графа
$X_{i+1}$ он содержит не более одного простого цикла.

Рассмотрим пару $(H_1,H_2)$, где граф $H_1$ является цепью длины
$\left[\frac{1}{\rho-1}\right]+1$, объединенной с графом,
изоморфным $\widetilde{X}_{i+1}^{r+1}$ (вершина $x_{i+1}^{i+1}$
при соответствующем изоморфизме переходит в некоторую вершину
$h$). Граф $H_2$ содержит лишь одну вершину, которая является
концевой вершиной рассмотренной цепи, отличной от $h$. Рассмотрим,
кроме того, множество $\mathcal{K}$ всех различных пар
$(K_1,K_2)$, для каждой из которых найдется такой граф $K$, что
$K\cap K_1=K_2$, граф $K\cup K_1$ является
$2^{k-1-i}$--расширением первого или второго типа графа $K$. В
силу $(n_1(\rho),n_2(\rho),n_3,n_4,\rho)$--разреженности графа
$Y_{i+1}$ найдутся такие вершины $y_{i+1}^{i+1}\in
V(Y_{i+1})\setminus
V(\widetilde{Y}_{i+1}^1\cup...\cup\widetilde{Y}_{i+1}^r)$,
$\widetilde{y}_{i+1}^l\in V(\widetilde{Y}_{i+1}^l)$ для некоторого
$l\in\{1,...,r\}$, цепь $c_{Y_{i+1}}\subset Y_{i+1}$ длины
$\left[\frac{1}{\rho-1}\right]+1$, соединяющая вершины
$\widetilde{y}_{i+1}^l$ и $y_{i+1}^{i+1}$, а также граф
$\widetilde{Y}_{i+1}^{r+1}\subset Y_{i+1}$, изоморфный
$\widetilde{X}_{i+1}^{r+1}$, что
$d_{Y_{i+1}}(y_{i+1}^{i+1},\widetilde{Y}_{i+1}^1\cup...\cup\widetilde{Y}_{i+1}^r)=
\left[\frac{1}{\rho-1}\right]+1$, $V(c_{Y_{i+1}})\cap
V(\widetilde{Y}_{i+1}^{r+1})=\{y_{i+1}^{i+1}\}$ и пара
$(c_{Y_{i+1}}\cup
\widetilde{Y}_{i+1}^{r+1},(\{\widetilde{y}_{i+1}^l\},\varnothing))$
является $(K_1,K_2)$--максимальным в $Y_{i+1}$ точным
$(H_1,H_2)$--расширением графа
$(\{\widetilde{y}_{i+1}^l\},\varnothing)$ для всех
$(K_1,K_2)\in\mathcal{K}$.\\

Продолжим изоморфизм графов $\varphi_{i+1}$ на вершины из
множества $V(\widetilde{X}_{i+1}^{r+1})$:
$\varphi_{i+1}|_{\widetilde{X}_{i+1}^{r+1}}:
\widetilde{X}_{i+1}^{r+1}\rightarrow \widetilde{Y}_{i+1}^{r+1}$,
причем $\varphi_{i+1}(x_{i+1}^{i+1})=y_{i+1}^{i+1}$.\\

Докажем, что графы
$\widetilde{X}_{i+1}^1,...,\widetilde{X}_{i+1}^{r+1},\widetilde{Y}_{i+1}^1,...,\widetilde{Y}_{i+1}^{r+1}$
обладают свойствами (I), (II) и (III) при $i<k-1$. Свойство (I)
выполнено, так как $x_{i+1}^{i+1}\in
V(\widetilde{X}_{i+1}^{r+1})$, $y_{i+1}^{i+1}\in
V(\widetilde{Y}_{i+1}^{r+1})$. Из изоморфности пар
$(c_{Y_{i+1}}\cup
\widetilde{Y}_{i+1}^{r+1},(\{\widetilde{y}_{i+1}^l\},\varnothing))$,
$(H_1,H_2)$ следует изоморфность пар
$(\widetilde{X}^1_{i+1}\cup...\cup
\widetilde{X}^{r+1}_{i+1},\widetilde{X}^1_{i+1}\cup...\cup
\widetilde{X}^{r}_{i+1})$, $(\widetilde{Y}^1_{i+1}\cup...\cup
\widetilde{Y}^{r+1}_{i+1},\widetilde{Y}^1_{i+1}\cup...\cup
\widetilde{Y}^{r}_{i+1})$ и справедливость свойства (III).
Осталось доказать, что графы
$\widetilde{X}^1_{i+1},...,\widetilde{X}^{r+1}_{i+1}$ обладают
$(k,i+1,r+1)$--свойством в $X_{i+1}$, а графы
$\widetilde{Y}^1_{i+1},...,\widetilde{Y}^{r+1}_{i+1}$ обладают
$(k,i+1,r+1)$--свойством в $Y_{i+1}$. Пусть $j\in\{1,...,r\}$.
Предположим, что существует такая вершина $v\in
V(\widetilde{X}_{i+1}^{r+1})$, отличная от $x_{i+1}^{i+1}$, что
$d_{X_{i+1}}(v,\widetilde{X}_{i+1}^j)<2^{k-1-i}$. Тогда в графе
$X_{i+1}$ существует подграф, являющийся $2^{k-i}$--расширением
первого типа графа $\widetilde{X}_{i+1}^j$. Получили противоречие
с $(k,i,r)$--свойством графов
$\widetilde{X}^1_{i+1},...,\widetilde{X}^{r}_{i+1}$ в $X_{i+1}$.
Следовательно,
$d_{X_{i+1}}(\widetilde{X}^j_{i+1},\widetilde{X}^{r+1}_{i+1})\geq
2^{k-1-i}$. Кроме того,
$$
 d_{Y_{i+1}}(\widetilde{Y}^j_{i+1},\widetilde{Y}^{r+1}_{i+1})>
 d_{Y_{i+1}}(\widetilde{Y}^j_{i+1},y_{i+1}^{i+1})-|V(\widetilde{Y}_{i+1}^{r+1})|\geq
$$
$$
 \geq\left[\frac{1}{\rho-1}\right]+1-2^{k-1-i}>2^{k-1}-2^{k-1-i}\geq 2^{k-2}\geq 2^{k-1-i}.
$$
Граф $\widetilde{X}^{r+1}_{i+1}$ (граф
$\widetilde{Y}_{i+1}^{r+1}$) не имеет общих вершин с графом
$\widetilde{X}_{i+1}^1\cup...\cup\widetilde{X}_{i+1}^r$ (графом
$\widetilde{Y}_{i+1}^1\cup...\cup\widetilde{Y}_{i+1}^r$) по
построению. Справедливы соотношения
$$
 |W_1\cup...\cup W_{r+1}|=
 |V(W_1\cup...\cup W_{r})|+|V(W_{r+1})|\leq
$$
$$
 \leq\eta(\rho)+(i-1)\left(\left[\frac{1}{\rho-1}\right]+1\right)+2^{k-1-i}<
 \eta(\rho)+i\left(\left[\frac{1}{\rho-1}\right]+1\right),
$$
где либо $W_j=\widetilde{X}^j_{i+1}$ для всех $j\in\{1,...,r+1\}$,
либо $W_j=\widetilde{Y}^j_{i+1}$ для всех $j\in\{1,...,r+1\}$.
Докажем, наконец, что в графе $X_{i+1}$ (в графе $Y_{i+1}$) не
найдется подграфа, являющегося $2^{k-1-i}$--расширением графа
$\widetilde{X}_{i+1}^{r+1}$ (графа $\widetilde{Y}_{i+1}^{r+1}$). В
случае графа $\widetilde{X}_{i+1}^{r+1}$ достаточно вспомнить, что
он содержит наибольшее количество вершин среди всех графов,
изоморфных какому-либо графу из $\mathcal{H}_{2^{k-1-i}}$ и
содержащих вершину $x_{i+1}^{i+1}$. Граф $Y_{i+1}$ не содержит
графов, являющихся $2^{k-1-i}$--расширениями первого или второго
типа графа $\widetilde{Y}_{i+1}^{r+1}$, так как граф
$\widetilde{Y}_{i+1}^{r+1}$ является $(K_1,K_2)$--максимальным в
$Y_{i+1}$ для всех $(K_1,K_2)\in\mathcal{K}$.\\

Если $i=k-1$, то так как $\varphi_k:
\widetilde{X}_k^1\cup...\cup\widetilde{X}_k^{r+1}\rightarrow
\widetilde{Y}_k^1\cup...\cup\widetilde{Y}_k^{r+1}$ --- изоморфизм
и $\varphi_k(x_k^j)=y_k^j$ при всех $j\in\{1,...,k\}$, то графы
$X_k|_{\{x_k^1,...,x_k^k\}}$, $Y_k|_{\{y_k^1,...,y_k^k\}}$ также
изоморфны и Консерватор побеждает.

\end{enumerate}

\section{Доказательство теоремы 1}
\label{proof1}

Пусть $\beta\in(0,\infty)\setminus \mathcal{Q}$ --- произвольное
число, $\alpha=1-\frac{1}{2^{k-1}+\beta}$, $\rho=1/\alpha$. Для
того, чтобы с вероятностью, стремящейся к 1, Консерватор выиграл в
игре $\EHR(G(N,N^{-\alpha}),G(M,M^{-\alpha}),k)$, пользуясь
стратегией, описанной в предыдущем разделе, достаточно доказать,
что с вероятностью, стремящейся к 1, случайный граф
$G(N,N^{-\alpha})$ является $(n_1,n_2,n_3,n_4,\rho)$--разреженным,
где числа $n_1,n_2,n_3,n_4$ определены равенствами (\ref{n1}),
(\ref{n4}). Обратимся сначала к первому свойству в определении
$(n_1,n_2,n_3,n_4,\rho)$--разреженного графа. Рассмотрим такое
множество $\mathcal{G}$ попарно неизоморфных графов, количество
вершин которых не превосходит $n_1$, а плотность отлична от
$\rho$, что любой граф $G$ с $v(G)\leq n_1$ и $\rho(G)\neq\rho$
изоморфен некоторому графу из $\mathcal{G}$. Пусть $\mathcal{G}_1$
--- такое множество попарно неизоморфных графов, количество вершин которых не превосходит
$n_1$, а максимальная плотность меньше, чем $\rho$, что любой
граф, удовлетворяющий заданным условиям, изоморфен некоторому
графу из $\mathcal{G}_1$. Очевидно, $|\mathcal{G}_1|\leq
|\mathcal{G}|<\infty$. Поэтому в соответствии с теоремой
\ref{erdos} справедливы равенства
$$
 \lim_{N\rightarrow\infty}{\sf P}_{N,p}(\forall G\in\mathcal{G}_1 \,\,
 N_G>0)=1,
$$
$$
 \lim_{N\rightarrow\infty}{\sf P}_{N,p}(\exists G\in\mathcal{G}\setminus\mathcal{G}_1 \,\,
 N_G>0)=0.
$$
Свойство 1) доказано.

Пусть $\mathcal{H}$ --- множество таких попарно неизоморфных
$1/\rho$--надежных пар $(H_1,H_2)$, что $v(H_1)\leq n_1$,
$v(H_2)\leq n_2$ и мощность $\mathcal{H}$ максимальна. Пусть
$\mathcal{K}$
--- множество таких попарно неизоморфных $1/\rho$--жестких пар $(K_1,K_2)$, что
$v(K_1)\leq n_3$, $v(K_2)\leq n_4$ и мощность $\mathcal{K}$
максимальна. В силу теоремы \ref{maximal_extensions_rigid} для
любых вершин $\widetilde{x}_1,...,\widetilde{x}_{v(H_2)}\in V_N$
выполнено равенство
$$
 \lim_{N\rightarrow\infty}{\sf P}_{N,p}(\forall (H_1,H_2)\in\mathcal{H} \, \forall (K_1,K_2)\in\mathcal{K}
 \,\,\,
 N_{(H_1,H_2)}^{(K_1,K_2)}(\widetilde{x}_1,...,\widetilde{x}_{v(H_2)})>0)=1.
$$
Тем самым, свойство 2) а, следовательно, и теорема \ref{positive}
доказаны.

\section{Доказательство теоремы 2}
\label{proof2}

В этом разделе для каждого $\alpha=1-\frac{1}{2^{k-1}+\beta},$
$\beta\in\widetilde{Q}$, мы рассмотрим два графа $G$ и $H$, а
также свойство $L_1$, которым обладает граф $G$, и свойство $L_2$,
которым обладает граф $H$. Мы докажем, что у Новатора есть
выигрышная стратегия в игре $\EHR(G,H,k)$, а затем установим, что
случайный граф $G(N,p)$ с вероятностью, стремящейся к некоторому
положительному числу, обладает свойствами $L_1$ и $L_2$. Тем самым
утверждение теоремы будет доказано.

Итак, пусть $k>3$ --- некоторое натуральное число,
$\beta\in\widetilde{Q}$, $\alpha=1-\frac{1}{2^{k-1}+\beta}$.
Положим
$\rho=\frac{1}{\alpha}=\frac{2^{k-1}+\beta}{2^{k-1}+\beta-1}.$
Рассмотрим граф $X$, равный объединению двух простых циклов
$C_1,C_2$, имеющие ровно одну общую вершину, причем
$v(C_1)+v(C_2)=2^{k-1}+\beta-1$, $v(C_1)\leq 2^{k-1}$, $v(C_2)\leq
2^{k-1}-1$, $v(C_1)>v(C_2)$. Пусть графы $C_1$ и $C_2$ не имеют
общих ребер. Обозначим $L_1$ свойство графа содержать подграф,
изоморфный $X$. Рассмотрим свойство содержать такую вершину $x$,
что в графе существует $v(C_1)$--расширение первого типа подграфа
$\{\{x\},\varnothing\}$ и $v(C_2)$--расширение первого типа того
же подграфа. Обозначим отрицание этого свойства $L_2$. Пусть граф
$G$ обладает свойством $L_1$, а граф $H$ обладает свойством $L_2$.
Опишем выигрышную стратегию Новатора в игре $\EHR(G,H,k)$. В
первом раунде Новатор выбирает такую вершину $x_1$ графа $G$, что
$x_1$ является вершиной некоторого подграфа $G_1\cup G_2$ в $G$,
изоморфного $X$, причем $G_1\cong C_1$, $G_2\cong C_2$,
$V(G_1)\cap V(G_2)=\{x\}$. Консерватор выбирает некоторую вершину
$y_1\in V(H)$. Не ограничивая общности, будем считать, что в графе
$H$ не найдется $v(C_1)$--расширения первого типа графа
$\{\{y_1\},\varnothing\}$. Новатор во втором раунде выбирает
вершину $x_2$ графа $G_1$, наиболее удаленную от вершины $x_1$.
Пусть $A$, $B$ --- две различные цепи, являющиеся подграфами $G_1$
и соединяющие $x_1$ с $x_2$. Консерватор во втором раунде выбирает
некоторую вершину $y_2\in V(H)$. Заметим, что в графе $H$ не
найдется двух различных цепей, длина одной из которых не
превосходит $v(A)-1$, а второй
--- $v(B)-1$, и соединяющих $y_1$ с $y_2$. Пусть, например, не
нашлось цепи длины $v(A)-1$. Положим $T_{1}^2=A$. Пусть сыграно
$i$ раундов, $i\geq 2$. Пусть, кроме того, выбраны вершины
$x_1,...,x_i\in V(G)$, $y_1,...,y_i\in V(H)$, а также такая цепь
$T_{i-1}^i$ соединяющая вершины $x_{i-1}$ и $x_i$ в графе $G$, что
в графе $H$ не найдется цепи такой же длины, соединяющей вершины
$y_{i-1}$ и $y_i$. Новатор в $i+1$-ом раунде выбирает такую
вершину $x_{i+1}$, принадлежащую цепи $T_{i-1}^i$, что
$|d_{T_{i-1}^i}(x_i,x_{i+1})-d_{T_{i-1}^i}(x_{i-1},x_{i+1})|$ ---
минимально. Консерватор во втором раунде выбирает некоторую
вершину $y_{i+1}\in V(H)$. Заметим, что в графе $H$ не найдется
двух цепей, одна из которых имеет длину
$d_{T_{i-1}^i}(x_{i-1},x_{i+1})$ и соединяет $y_{i-1}$ с
$y_{i+1}$, вторая имеет длину $d_{T_{i-1}^i}(x_{i},x_{i+1})$ и
соединяет $y_{i}$ с $y_{i+1}$ (мы обозначаем $d_W(x,y)$ длину
наименьшей цепи в графе $W$, соединяющей вершину $x$ с вершиной
$y$). Пусть, например, не нашлось цепи, соединяющей $y_{i-1}$ с
$y_{i+1}$. Обозначим цепь, являющуюся подграфом в $T_{i-1}^i$ и
соединяющую $x_{i-1}$ с $x_{i+1}$, $T_i^{i+1}$. Так как
$\max\{v(C_1),v(C_2)\}\leq 2^{k-1}$, то в одном из раундов с
номером $r\in\{3,...,k\}$ Новатор выберет такую вершину $x_r$, что
она будет соединена ребрами как с вершиной $x_{r-1},$ так и с
вершиной $x_{r-2}$. В графе $H$ вершины, соединенной и с
$y_{r-1}$, и с $y_{r-2}$, не найдется. Поэтому Новатор победит.

Докажем теперь, что случайный граф $G(N,p)$ с вероятностью,
стремящейся к некоторому положительному числу, обладает свойствами
$L_1$ и $L_2$. Граф $X$ строго сбалансированный, а его плотность
равна $\rho$. Поэтому в силу теоремы \ref{erdos} вероятность того,
что случайный граф обладает свойством $L_1$ стремится к
$1-e^{-1/a(X)}$. Для завершения доказательства осталось установить
справедливость следующего соотношения:
$\lim_{N\rightarrow\infty}{\sf P}_{N,p}(\overline{L}_2)\in(0,1)$.
Легко заметить, что свойство $\overline{L}_2$ выполнено тогда и
только тогда, когда граф содержит подграф, выбранный из некоторого
конечного множества, причем плотность графов из этого множества не
меньше, чем $\rho$. В соответствии с теоремой о совместном
распределении чисел подграфов в случайном графе (см.
\cite{Janson}, Chapter III, Remark 3.20) с вероятностью,
стремящейся к некоторому числу из интервала $(0,1)$, случайный
граф $G(N,p)$ содержит хотя бы один подграф из упомянутого
множества. Тем самым теорема доказана.

\end{document}